\documentclass[11pt,reqno]{amsart}
\usepackage{amsfonts}
\usepackage{fancyhdr}
\usepackage[dvips]{graphics}
\usepackage{epsfig}
\usepackage{epsf}
\usepackage{hyperref}
\usepackage{flafter}
\usepackage{amsmath,amsfonts,amsthm,amssymb,amscd}
\usepackage{color}
\usepackage[parfill]{parskip}
\usepackage{epstopdf}

\usepackage{verbatim}

\newtheorem{thm}{Theorem}

\newtheorem{lem}[thm]{Lemma}
\newtheorem*{thm*}{Theorem}
\newtheorem{con*}{Conjecture}
\newtheorem*{lem*}{Lemma}
\newtheorem{prop}[thm]{Proposition}
\newtheorem{defn}{Definition}[section]
\newtheorem{rem*}{Remark}

\def\1{x_1}
\def\2{x_2}
\def\3{x_3}
\def\4{x_4}

\begin{document}
 \baselineskip=17pt
\hbox{}
\medskip

\title
[index - singular case] {On the index conjecture in zero-sum theory: \\ singular case}

\author{Fan Ge}

\email{fange.math@gmail.com}

\address{Department of Mathematics, University of Rochester, Rochester, NY, United States}

\address{Current Address: Department of Pure Mathematics, University of Waterloo, Waterloo, ON, Canada}

\begin{abstract}
Let $S=(a_1)\cdots(a_k)$ be a minimal zero-sum sequence over a
finite cyclic group $G$ of order $n$. An important question in zero-sum theory is to determine the pairs
$(k, n)$ for which every minimal zero-sum sequence $S$ of length $k$
over $G$ has index 1. Progress towards this question has been made by various authors; the only case that remains open, known as the index conjecture, is when $k=4$
and $\gcd(|G|,6)=1$.  In this paper we make a contribution to the index conjecture. Namely, we prove
that if $S$ is singular then the index of $S$ is $1$.
\end{abstract}

\maketitle

\section{Introduction}

Throughout this paper $G$ is a finite additive cyclic group of order $n$. By a \emph{sequence $S$ of length $k$ over $G$}
we mean a sequence with $k$ elements, each of which is in $G$. We
write $(a_1)\cdots(a_k)$ for such a sequence. A sequence $S$ is a \emph{zero-sum}
sequence if $\sum_i a_i=0$. If $S$ is a zero-sum sequence but no proper nontrivial subsequence of
it is, then we say $S$ is a \emph{minimal} zero-sum sequence. Given any generator $g$ of $G$, we may write
$S=(x_1g)\cdots(x_kg)$ for some natural numbers $x_1, \dots, x_k$,
where by $x_ig$ we mean the sum $g+g+ \cdots +g$ with $x_i$ terms.

\begin{defn}\textnormal{Let $S=(x_1g)\cdots(x_kg)$ be a sequence over $G$,
where $1\le x_1,...,x_k\le n$.
Define the $g$-norm of $S$ to be
$\|S\|_g=\frac{\sum_{i=1}^{k}x_i}{n}$. The \emph{index} of $S$ is
defined by
$$\textnormal{ind}(S)=\min\|S\|_g,$$ where the minimum is taken
over all generators $g$ of $G$.}
\end{defn}

The index of a sequence is an important invariant in zero-sum theory. It plays a crucial role in the
study of zero-sum sequences and related topics (see, for example, Geroldinger~\cite{g90}  and Gao~\cite{gao}).
An important question is to determine the pairs
$(k, n)$ for which every minimal zero-sum sequence $S$ of length $k$
over $G$ has index 1. The cases  $k\ne 4$ or
$\gcd(n, 6)\ne 1$ have been settled (see~\cite{Gero}, \cite{P}, \cite{sc}, \cite{Y}).
Therefore, the only remaining case   is when both  $k=4$ and $\gcd(n,
6)=1$. The following conjecture is widely held.

\begin{con*}\label{conj}
Let $G$ be a finite cyclic group such that $\gcd(|G|, 6) = 1$. Then
every minimal zero-sum sequence $S$ over $G$ of length $4$ has
\textnormal{ind}($S$) = $1$.
\end{con*}

\begin{rem*}
\textnormal{It is easy to see that, for such $S$  we have either $\textnormal{ind}(S)=1$ or
$\textnormal{ind}(S)=2$, and moreover, $\textnormal{ind}(S)=2$ if
and only if $\|S\|_g=2$ for all generators $g$ of $G$. Indeed, for such $S$ we have by definition that $\|S\|_g$ could be 1, 2, or 3; but if $\|S\|_g=3$ for some generator $g$, then $\|S\|_{-g}=1$, where $-g$ is also a generator.}
\end{rem*}

Below we always assume that $(n, 6)=1$.

In~\cite{LPYZ}, Y. Li et al. proved that if $n$ is a prime power then Conjecture~\ref{conj} is
true. Later the case
when $n$ has two distinct prime factors was also proved (see~\cite{LP}
and~\cite{XS}). Recently, X. Zeng and X. Qi \cite{ZQ} proved the conjecture for the case when $n$ is coprime to $30$. In~\cite{G} we proved the following result.

\begin{thm*}
Suppose that $n$ is the smallest integer for which
Conjecture~\ref{conj} fails. Let $S=(x_1)(x_2)(x_3)(x_4)$ be a
minimal zero-sum sequence over $G\cong \mathbb{Z}/n$ with
\textnormal{ind}$(S)=2$. Then we have $\gcd(n, x_i)=1$ for all $i$.
\end{thm*}

In view of the above theorem, we see that to prove Conjecture~\ref{conj} it suffices to prove that
if $S=(x_1)(x_2)(x_3)(x_4)$ is a minimal zero-sum sequence
 with  $(n, x_i)=1$ for all $i$, then \textnormal{ind}$(S)=1$.

The purpose of this paper is to give a proof of
Conjecture~\ref{conj} in the case when $S$ is singular (see below).

\begin{defn}\textnormal{Let $S=(x_1)\cdots(x_4)$ be a minimal zero-sum sequence over $\mathbb{Z}/n$ with $1\le x_1,...,x_4< n$ and $(n, 6)=1$.
Then $S$ is called \emph{singular} if the following conditions hold:\\
(i) $x_1=1$; \\
(ii) $x_2+1=x_3$ or $x_2=n-2$; \\
(iii) $(x_i, n)=1$ for all $i=1,\dots,4$.
}
\end{defn}

\begin{thm}\label{thm singular 1}
The index of $S$ is $1$ if $S$ is singular.
\end{thm}

\section{preliminaries}

Throughout this section we always assume that $S=(x_1)\cdots(x_k)$ is a minimal zero-sum sequence over $\mathbb{Z}/n$ with $1\le x_1,...,x_k< n$.

For integers $x$ and $y>0$, let $(x)_y$ denote the least nonnegative
residue of $x$ mod $y$. For $z\in \mathbb Z/y$, we may view $z$ as
an integer and define $(z)_y$ similarly.

\begin{lem} \label{lem coeff} Let $S=(x_1)(x_2)(x_3)(x_4)$.
Given any generator $g$ in $G$, write $S=(y_1g)(y_2g)(y_3g)(y_4g)$
for $1\le y_i < n$. Then we have $y_i=(g^{-1}x_i)_n$ for $i=1,
\dots, 4$, where $g^{-1}$ is the inverse of $g$ in the multiplicative
group $(\mathbb Z/n)^*$. In particular, we have
$$
n\ \|S\|_{g^{-1}}\ =\ \sum_{i=1}^4 (gx_i)_n\ .
$$
\end{lem}

\proof For any $i=1, ..., 4$, we have $x_i=y_i g$ in $G$. Hence
$x_ig^{-1}=y_i$ in $G$. Therefore, $y_i=(g^{-1}x_i)_n$.
\qed

The following lemma will be used frequently.

\begin{lem}\label{lem 0}
If $S$ has index $2$, then \\
(i) For any $g\in (\mathbb Z/n)^*$ we have
$
\#\{i: (x_ig)_n>\frac{n}{2}\}=2.
$
\\
(ii) If in addition $(x_i, n)=1$ for all $i$, then the $x_i$'s are pairwise distinct.
\end{lem}
\proof The first part is essentially Remark 2.1 of \cite{LP}. For the second part,
suppose $x_i=x_j$ for some $i\ne j$. By Remark 1 and Lemma~\ref{lem coeff} we have
$$
2n=n\ \|S\|_{x_i}\ = \sum_{u=1}^4 (x_i^{-1}x_u)_n= 1+1+\sum_{u\ne i,j}(x_i^{-1}x_u)_n.
$$
Since $S$ is minimal, we clearly have $(x_i^{-1}x_u)_n\le n-2$. Thus, the right-hand side of the above is at most $1+1+(n-2)+(n-2)=2n-2$, a contradiction.
\qed

To state our next result, we make the following definition.

\begin{defn}\textnormal{We call an integer $k$ \emph{good} if $k$ satisfies the following conditions:\\
(i) $k=2^l$, where $l$ is a nonnegative integer;\\
(ii) $k<\frac{n}{6}$;\\
(iii) $F(k):=\Big( 2n-2-2\big[ \frac{3k-1}{3k}\,n \big]\Big)\,k>\frac{n-1}{2}$, where $[\cdot]$ is the floor function.
}
\end{defn}

\begin{prop} \label{prop induction}
Let $S=(x_1)(x_2)(x_3)(x_4)$. Suppose that $S$ has index $2$, and that $x_1=1$, $x_2+1=x_3$ and $(x_i, n)=1$ for all $i$.
If $k$ is good, then $x_2\ge [\frac{6k-1}{6k}n]$.
\end{prop}

We need the following three lemmas in order to prove Proposition \ref{prop induction}.

\begin{lem}\label{lem x4<x2}
Under the assumption of Proposition \ref{prop induction}, we have $\4<\2$.
\end{lem}

\proof

Since $\1+\4=1+\4<n$ we see that $2\2+1=\2+\3\ge n+1$. It follows that $\2\ge \frac{n+1}{2}$.
By part (i) of Lemma~\ref{lem 0} we have
$
\#\{i: x_i>\frac{n}{2}\}=2.
$
Therefore, as $\3>\2>n/2$, we must have $\4<n/2$. Hence $\4<\2$.
\qed

\begin{lem} \label{lem good}
If $k\ge 2$ is good, then $k/2$ is good.
\end{lem}

\proof

It suffices to check the third condition. Namely, we need to show that
$$
F\Big(\frac{k}{2}\Big)>\frac{n-1}{2}.
$$
Since $k$ is good, we have
$$F(k)>\frac{n-1}{2}.$$
Therefore it suffices to prove that $F(k/2)\ge F(k)$.
Writing
\begin{align}\label{f}
f(k)=\Big[ \frac{3k-1}{3k}\,n \Big],
\end{align} then by a straightforward computation, we see that to prove $F(k/2)\ge F(k)$
it suffices to show that
$$
2f(k)-f(k/2) \ge n-1.
$$
Notice that the left-hand side is an integer. Thus, the above inequality is equivalent to
$$
2f(k)-f(k/2) > n-2,
$$
or,
$$
2\Big(\frac{3k-1}{3k}n -\Big\{\frac{3k-1}{3k}n\Big\}\Big)-\Big(\frac{3k/2-1}{3k/2}n -\Big\{\frac{3k/2-1}{3k/2}n\Big\}\Big) > n-2.
$$
A straightforward computation turns it into
$$
2-2\Big\{\frac{3k-1}{3k}n\Big\}+\Big\{\frac{3k/2-1}{3k/2}n\Big\}>0.
$$
But this is clearly true. Hence the lemma follows.
\qed

\begin{lem} \label{lem induct}
Let $f(k)$ be defined as in~\eqref{f}. Under the assumption of Proposition \ref{prop induction}, if $k$ is good and $x_2\ge f(k)$, then $x_2\ge f(2k)$.
\end{lem}

\proof

First we claim that $x_2\ne f(k)$.
Suppose $x_2=f(k)$. We will show that
$$
\#\Big\{i: (x_ik)_n>\frac{n}{2}\Big\}\ne 2.
$$
Then by Lemma~\ref{lem 0} $S$ would have index $1$, a contradiction.
To prove the above inequality, first notice that $\1 k=k$. Thus $(\1 k)_n=k$.
Next,
$$
\2 k=f(k)k=\Big[\frac{3k-1}{3k}n\Big]k\in \Big(\,\Big(\frac{3k-1}{3k}n-1\Big)k, \frac{3k-1}{3k}nk\Big),
$$
and
$$
\Big(\frac{3k-1}{3k}n-1\Big)k-(k-1)n=2n/3-k>0
$$
since $k$ is good implies $k<n/6$. Hence, $\2k\in((k-1)n,kn)$, and this gives
\begin{align*}
(\2k)_n=\2k-(k-1)n&>\Big(\frac{3k-1}{3k}n-1\Big)k-(k-1)n=2n/3-k\\
& >2n/3-n/6=n/2.
\end{align*}
Also, since $\3=\2+1$, we have $\3k\in(\2k,kn)\subseteq ((k-1)n,kn)$. Thus, $$(\3k)_n=\3k-(k-1)n>x_2k-(k-1)n=(\2k)_n>n/2.$$
Next, since $\4=2n-\1-\2-\3=2n-2-2\2=2n-2-2f(k)$, and since $k$ is good implies that $(2n-2-2f(k))k>\frac{n-1}{2}$, we see that $\4k>\frac{n-1}{2}$.
It follows that $\4k\ge \frac{n+1}{2}$ since $\4k$ is an integer. On the other hand,
$$
\4k=(2n-2-2f(k))k<\Big(2n-2-2\Big(\frac{3k-1}{3k}n-1\Big)\Big)k=2n/3.
$$
Therefore, we see that $(\4k)_n=\4k>n/2$. This gives
$$
\#\Big\{i: (x_ik)_n>\frac{n}{2}\Big\}=3,
$$
and the claim follows.

Thus, we have $\2\ge f(k)+1$. It then follows that
\begin{align*}
\2 3k & \in ((f(k)+1)3k,3kn)\\
& = \Big(\,\Big(\Big[\frac{3k-1}{3k}n\Big]+1\Big)3k, 3kn\Big)\\
& \subseteq\Big(\frac{3k-1}{3k}n\cdot 3k, 3kn\Big)\\
&=((3k-1)n, 3kn).
\end{align*}
This implies that $(\23k)_n=\23k-(3k-1)n$.
We also clearly have
$$
\33k\in(\23k,3kn)\subseteq((3k-1)n, 3kn).
$$
Hence, $(\33k)_n=\33k-(3k-1)n$.
Recall that $k<n/6$, thus $(\13k)_n=3k$.
Now since $S$ has index $2$ and since $\gcd(3k, n)=1$, by Remark 1 and Lemma~\ref{lem coeff} we see that
$$
\sum_{i=1}^4(x_i 3k)_n=2n.
$$
It follows that
\begin{align*}
2n&=3k+\23k-(3k-1)n+\33k-(3k-1)n+(\43k)_n\\
& = 3k(\1+\2+\3+\4)-(6k-2)n+(\43k)_n-\43k\\
& =3k\cdot 2n-(6k-2)n+(\43k)_n-\43k\\
& = 2n+(\43k)_n-\43k.
\end{align*}
Therefore, we have $\43k=(\43k)_n<n$, or $\4<n/3k$. Using the relation $\4=2n-2-2\2$, we obtain
$$
\2>n-1-\frac{n}{6k}.
$$
But since $\frac{n}{6k}$ is not an integer, we conclude that $\2\ge[n-n/6k]=f(2k)$.
This completes the proof of the lemma.
\qed

\proof[Proof of Proposition~\ref{prop induction}:]
First we show that $\2\ge f(1)$. By Lemma~\ref{lem x4<x2} we have $\4<\2$. But $\4=2n-2-2\2$. It follows that
$2n-2-2\2<\2$, or $\2>\frac{2n-2}{3}$. When $n\equiv 1\pmod 3$, say $n=3m+1$, we have $\2>\frac{2n-2}{3}=2m=[2n/3]=f(1)$.
When $n\equiv 2\pmod 3$, say $n=3m+2$, we have $$\2>\frac{2n-2}{3}=2m+\frac{2}{3}.$$ This implies $\2\ge 2m+1$ since $\2$ is an integer.
But $2m+1=[2n/3]=f(1)$. Thus, in both cases we have $\2\ge f(1)$.

Write $k=2^l$. Since $k$ is good, by Lemma~\ref{lem good} we have $2^t$ is good for any integer $t\in [0,l]$.
In particular, $2^0=1$ is good. This together with the fact that $\2\ge f(1)$ implies $\2\ge f(2)$ by Lemma~\ref{lem induct}.
Now since $\2\ge f(2)$ and since $2^t$ is good for any integer $t\in [0,l]$, we can use Lemma~\ref{lem induct} repeatedly to conclude that
$\2\ge f(2\cdot 2^l)=f(2k)$.
\qed

\begin{lem}\label{lem k value}
Let $b\ge 3$ be an integer such that $$3\cdot 2^b<n<3\cdot 2^{b+1}.$$ Then $k=2^{b-2}$ is good.
\end{lem}

\proof

Clearly $k=2^{b-2}$ satisfies the first two conditions in the definition of good. It remains to prove that
$$\Big( 2n-2-2\big[ \frac{3k-1}{3k}\,n \big]\Big)\,k>\frac{n-1}{2}.$$
It suffices to show that
$$\Big( 2n-2-2\Big(\frac{3k-1}{3k}\,n \Big)\Big)\,k>\frac{n-1}{2}.$$
A straightforward computation shows that this is equivalent to
$$
2k<\frac{n}{6}+\frac{1}{2}.
$$
But this is clear since $$2k=2^{b-1}=\frac{3\cdot 2^b}{6}<\frac{n}{6}.$$
\qed

\begin{prop} \label{prop x2>=n-4}
Let $S=(x_1)(x_2)(x_3)(x_4)$. Suppose that $S$ has index $2$, and that $x_1=1$, $x_2+1=x_3$ and $(x_i, n)=1$ for all $i$. Then $x_2=n-4$ or $n-3$.
Therefore, $$S=(1)(n-4)(n-3)(6) \quad or \quad (1)(n-3)(n-2)(4).$$
\end{prop}

\proof
Let $3\cdot 2^b<n<3\cdot 2^{b+1}.$ By Proposition~\ref{prop induction}, we have $\2\ge f(2k)$ if $k$ is good. By Lemma~\ref{lem k value}, $k=2^{b-2}$ is good.
It follows that $$\2\ge f(2^{b-1})=\Big[\frac{3\cdot 2^{b-1}-1}{3\cdot 2^{b-1}}n\Big].$$
This implies that
\begin{align*}
n-\2 & < n-\frac{3\cdot 2^{b-1}-1}{3\cdot 2^{b-1}}n+1
\\
& < 5.
\end{align*}
Hence, $n-\2\le 4$, and $\2$ could be $n-4$, $n-3$ or $n-2$. The case $\2=n-2$ can be excluded because $n-2\ge \3=\2+1>\2$.
\qed

\section{Proof of Theorem~\ref{thm singular 1}}

In this section we prove Theorem~\ref{thm singular 1}. We need the following result.

\begin{thm} \label{thm explicit case index=1}
The sequences $(1)(n-4)(n-3)(6)$ and $(1)(n-3)(n-2)(4)$ have index $1$.
\end{thm}

\proof

If $n\le1000$, it is known that every length four minimal zero-sum sequence over $\mathbb Z/n$ has index $1$ (see~\cite{P}, and also~\cite{LP}). Thus, we may assume $n> 1000$.
 Moreover, it is known that $S$ has index $1$ if $n$ has at most two distinct prime factors (see~\cite{LPYZ}, ~\cite{LP} and~\cite{XS}).
 Hence we may also assume that $n$ has at least three distinct prime factors.

First, let us consider $S=(1)(n-4)(n-3)(6)$.
We claim that there exists a $g\in (\mathbb Z/n)^*$ such that $n/12<g<n/8$. This is clearly true for $n$ ``large enough'' in view of the
Prime Number Theorem in arithmetic progressions. But here we would like to avoid determining the effective lower bound for ``large'' $n$.
Therefore, instead of the Prime Number Theorem, we use the following result.

\begin{lem}\label{lem 2n 3n}
\textnormal{(See \cite{B}.)}\\
\textnormal{(i)} The interval $(2N, 3N)$ contains a prime for any integer $N\ge 2$.\\
\textnormal{(ii)} The interval $\big[N+1, \frac{3(N+1)}{2}\big)$ contains a prime for any integer $N\ge 2$.
\end{lem}

If there is a prime $p$ in the interval $(n/12, n/8)$, then we will be done if we can show that $p$ is coprime to $n$. Indeed,
if $p$ divides $n$, then we would have $n/p$ is an integer less than $12$, and thus $n$ takes the form
$5p$, $7p$ or $11p$, which contradicts our assumption that $n$ has at least three distinct prime factors.

Now suppose there is no prime in the interval $(n/12, n/8)$.
Then in this case we can show that $[n/12]$ is a prime. Indeed, taking $N=[n/24]$ in part (i) of the above lemma, we see that the interval
$(2[n/24],3[n/24])$ contains a prime. It follows that $([n/12]-1, n/8)$ contains a prime since $2[n/24]\ge [n/12]-1$ and $3[n/24]\le n/8$.
Thus, if $(n/12, n/8)$ contains no prime then $[n/12]$ must be prime.

Let $[n/12]=q$. Since $q$ is prime and $n>1000$ we see that $q$ is odd.
In part (ii) of the above lemma take $N=q$, and it follows that $[q+1, 3(q+1)/2)$ contains a prime. Since $q$ is odd, $3(q+1)/2$ is an integer.
Thus, we obtain that the interval
$$\Big[q+1, \frac{3(q+1)}{2}-1\Big]$$ contains a prime. Note that $q+1=[n/12]+1>n/12$.
Therefore, using the fact that $[q+1, 3(q+1)/2-1]$ contains a prime and the assumption that $(n/12, n/8)$ contains no prime,
we conclude that $$\frac{3(q+1)}{2}-1 \ge \frac{n}{8}.$$
This in turn gives $12q+4\ge n$. Recall that $q=[n/12]$ and $\gcd(n,6)=1$. Hence we have $$n=12q+1.$$
Moreover, it is straightforward to compute that $$\frac{3(q+1)}{2}-2 - \frac{n}{8}=-\frac{5}{8}.$$
Thus we have $[q+1, 3(q+1)/2-2]\subseteq (n/12,n/8)$. It follows that $3(q+1)/2-1$ is a prime.

Write $q=2z+1$. Then since $3z+2=3(q+1)/2-1$ is prime we see that $z$ is odd. By part (ii) of the above lemma, we know that
$\big[z+1, \frac{3(z+1)}{2}\big)$ contains a prime. Hence $$\Big[z+1, \frac{3(z+1)}{2}-1\Big]$$ contains a prime $p$ since $3(z+1)/2$ is an integer.
We easily check that $$z+1>\frac{24z+13}{24}=\frac{n}{24}$$ and that $$\frac{3(z+1)}{2}-1<\frac{24z+13}{16}=\frac{n}{16}.$$
Thus, we conclude that $$p\in \Big(\frac{n}{24}, \frac{n}{16}\Big).$$
This gives $$2p\in \Big(\frac{n}{12}, \frac{n}{8}\Big).$$

To prove our claim, it remains to show that $\gcd(2p,n)=1$. If not, then $p$ divides $n$. But this implies that $n/p$ is an integer smaller than $24$, and thus $n$ takes the form
$5p, 7p, 11p, 13p, 17p, 19p$ or $23p$, which contradicts our assumption that $n$ has at least three distinct prime factors.
Therefore, our claim follows.

We have shown that there exists a $g\in (\mathbb Z/n)^*$ with $n/12<g<n/8$. A straightforward computation then shows that
$(1\cdot g)_n<n/2$, $((n-4)\cdot g)_n>n/2$, $((n-3)\cdot g)_n>n/2$ and that $(6\cdot g)_n>n/2$. Thus, we have
$$
\#\Big\{i: (x_ig)_n>\frac{n}{2}\Big\}=3,
$$
where $\1=1$, $\2=n-4$, $\3=n-3$ and $\4=6$. It follows from Lemma~\ref{lem 0} that $S$ has index $1$.

The proof for the case when $S=(1)(n-3)(n-2)(4)$ is almost exactly the same, and we omit it.
\qed

\medskip

\proof[Proof of Theorem~\ref{thm singular 1}:]

If $\3=\2+1$ and $S$ has index $2$, then by Proposition~\ref{prop x2>=n-4}
$$S=(1)(n-4)(n-3)(6) \quad or \quad (1)(n-3)(n-2)(4).$$
But it follows from Theorem~\ref{thm explicit case index=1} that these $S$ both have index $1$, a contradiction.
Thus, if $\3=\2+1$ then $S$ has index $1$.

Now suppose $\2=n-2$ and $S$ has index $2$. Then $2n=1+(n-2)+\3+\4$ implies that $(\3^{-1}\4)_n=(\3^{-1}-1)_n$.
Thus, we have
\begin{align*}
S &= (1)(\2)(\3)(\4)\\
& = (\3^{-1}\3)(\3^{-1}\2\3)(\3)(\3^{-1}\4\3)\\
& =\Big((\3^{-1})_n\3\Big)\Big((\3^{-1}\2)_n\3\Big)\Big((1)_n\3\Big)\Big((\3^{-1}\4)_n\3\Big)\\
&=\Big((\3^{-1})_n\3\Big)\Big((\3^{-1}\2)_n\3\Big)\Big((1)_n\3\Big)\Big((\3^{-1}-1)_n\3\Big)\\
&=:(y_3\3)(y_4\3)(y_1\3)(y_2\3),
\end{align*}
where $y_1=1=(1)_n$, $y_2=(\3^{-1}-1)_n$, $y_3=(\3^{-1})_n$ and $y_4=(\3^{-1}\2)_n$. Let $Y$ be the sequence $(y_1)(y_2)(y_3)(y_4)$.
It is clear that $Y$ is singular, and that $y_2+1=y_3$. Therefore, it follows from the last paragraph that $Y$ has index $1$.
It is then clear that $S$ also has index $1$.
This completes the proof.
\qed

\section*{acknowledgement}
This work was partially supported by NSF grant DMS-1200582.

\end{document}